\documentclass[journal]{ieeetran}
\usepackage{cite}
\usepackage{graphicx}
\usepackage{epstopdf}
\usepackage{amsmath}
\usepackage{cases}
\usepackage[caption=false,font=footnotesize]{subfig}
\usepackage{fixltx2e}
\usepackage{amssymb}
\usepackage{mathtools}
\usepackage{bm}
\usepackage{multirow}
\usepackage{color}
\usepackage[dvipsnames]{xcolor}
\usepackage{algorithm}
\usepackage{algorithmic}
\usepackage{amsthm}
\usepackage{mathrsfs}
\usepackage{textcomp,mathcomp}
\usepackage{booktabs}
\usepackage{threeparttable}

\newtheorem{remark}{Remark}
\usepackage{pifont}
\usepackage{threeparttable}
\usepackage[ruled,vlined,linesnumbered,algo2e]{algorithm2e}

\newcommand{\tabincell}[2]{\renewcommand\arraystretch{0.9}\begin{tabular}{@{}#1@{}}#2\end{tabular}}

\begin{document}

\title{Beyond Higher-Pulse Rectification: Operational Harmonic Coordination in Renewable P2H Systems}

\author{
Yangjun~Zeng,~\IEEEmembership{Student Member,~IEEE},
Yiwei~Qiu,~\IEEEmembership{Member,~IEEE},
Jie~Zhu,~\IEEEmembership{Student Member,~IEEE},
Yi~Zhou,~\IEEEmembership{Member,~IEEE},
Tao~Wu,~\IEEEmembership{Member,~IEEE},
Xin~Meng,~\IEEEmembership{Member,~IEEE},
Shi~Chen,~\IEEEmembership{Member,~IEEE},\\
Buxiang~Zhou,~\IEEEmembership{Member,~IEEE},
and
Kaigui~Xie,~\IEEEmembership{Fellow,~IEEE}%

\thanks{This work was supported in part by the National Natural Science Foundation of China under Grants 52377116 and 52577129, and in part by the Doctoral Student Program of the Young S\&T Talents Cultivation Project, CAST. All authors are with the College of Electrical Engineering, Sichuan University, Chengdu 610065, China. \emph{(Corresponding author: Yiwei Qiu) (email: ywqiu@scu.edu.cn)}}
%\thanks{The authors are with the College of Electrical Engineering, Sichuan University, Chengdu 610065, China. (email: ywqiu@scu.edu.cn)}
%
}
\maketitle

\begin{abstract}
	Thyristor rectifiers (TRs) are cost-effective electrolysis power supplies for renewable power-to-hydrogen (ReP2H) systems, but their harmonics may violate grid-code limits. In contrast to conventional solutions that rely on higher-pulse (such as 24-pulse) rectifiers, this paper proposes an operational harmonic coordination scheme that enables low-cost 12-pulse TRs to meet harmonic requirements through coordinated operation. First, a harmonic model quantifies the effects of rectifier transformer (RCT) tap positions and electrolytic currents, enabling harmonic cancellation among multiple electrolyzers (ELZs). A two-layer framework then coordinates hydrogen production and harmonic mitigation.
	Hourly scheduling determines ELZ commitment within the harmonic feasible region under renewable uncertainty using stochastic programming and a modified progressive hedging algorithm, while minute-level dispatch tracks renewable power and mitigates harmonics. A decomposition algorithm separates production dispatch from harmonic mitigation to improve computational efficiency. Case studies based on real-life projects show that the proposed method increases profit by 31\% over current-only regulation. Annual simulations further show that coordinated 12-pulse TRs can replace 24-pulse rectifiers for harmonic compliance by exchanging additional RCT tap actions for lower transformer investment, reducing rectification-stage cost by 37.5\%.

\end{abstract}


\section*{Nomenclature}
{\subsection{Abbreviations}
	\begin{IEEEdescription}[\IEEEusemathlabelsep\IEEEsetlabelwidth{superscript}]
		\addcontentsline{toc}{section}{Nomenclature}
		\item[ELZ] Electrolyzer
		\item[NWT] Network transformer
		\item[OLTC] On-load tap changer
		\item[PV] Photovoltaic
		\item[RCT] Rectifier transformer
		\item[ReP2H] Renewable power-to-hydrogen
		%\item[SP] Stochastic programming
		\item[TR] Thyristor rectifier
		\item[WT] Wind turbine
\end{IEEEdescription}}

\vspace{-8pt}
\subsection{Indices and Sets}
\begin{IEEEdescription}[\IEEEusemathlabelsep\IEEEsetlabelwidth{superscript}]
	\addcontentsline{toc}{section}{Nomenclature}
	
	%\item[$d$] Index for typical day in the planning framework
	%\item[$h$] Index for scenarios of renewable power and hydrogen price in each typical day
	%\item[$t$] Index for hourly scheduling time periods
	\item[$t$, $\tau$] Index for hourly and minute dispatch periods
	\item[$n,h$] Index for ELZs and harmonic order
	\item[$s$] Index for representative scenarios
	%\item[$m$]  Index for iteration of the progressive hedging
	\item[$j, j',j''$] Index for buses
	\item[$j'j$]   Index for branches
	%\item[$i$]  Index for all subproblem
	%\item[$\pi(j)$, $\sigma(j)$] Set of parents and children of bus $j$
\end{IEEEdescription}

\vspace{-6pt}
\subsection{Variables}
\subsubsection{ELZ-Related Operational Variables}
\begin{IEEEdescription}[\IEEEusemathlabelsep\IEEEsetlabelwidth{superscript}]
	\addcontentsline{toc}{section}{Nomenclature}
	\item[$P^{\text{Stack}}, U^{\text{Stack}}$] Electrolytic power and voltage
	\item[$I, Y^{\text{H}_2}$] Electrolytic current and hydrogen output rate
	\item[$T$] Stack temperature of ELZ
	\item[$U^{\text{AC}}, I^{\text{st}}$] AC-side voltage and fundamental current
	\item[${\textbf{I}}_h, \hat{\textbf{I}}_h$] $h$th harmonic current and its component
	\item[$K^{\text{rec}}, k^{\text{rec}}$] RCT turns ratio and OLTC tap position
	\item[$b^{\text{On}}, b^{\text{By}}, b^{\text{Idle}}$] Production, standby, and idle states of ELZ
	\item[$b^{\text{SU}}, b^{\text{SD}}$] Startup/shutdown actions of ELZ
	\item[$P^{\text{ELZ}}, P^{\text{BoP}}$] Active power of ELZ/balance of plant
	\item[$P^{\text{Loss}}$] The active loss of the TR of ELZ
	%\item[$P^{\text{Rec}}$] Active power of the rectifier of the ELZ
	\item[$Q^{\text{ELZ}}$] Reactive power of the TR in ELZ
	\item[$P^{\text{Heat}}, P^{\text{Diss}}$] Electrolytic heat and dissipation of ELZ
	\item[$P^{\text{Cool}}$] Cooling heat flow of ELZ
\end{IEEEdescription}

\subsubsection{Network-Side Operational Variables}
\begin{IEEEdescription}[\IEEEusemathlabelsep\IEEEsetlabelwidth{superscript}]
	\addcontentsline{toc}{section}{Nomenclature}
	%\item[$P_{ij,t}, Q_{ij,t}$] Active/reactive power flows on branch $ij$
	%\item[$I_{ij,t}, \ell_{ij,t}$] Current on branch $ij$ and its square
	%\item[$U_{j,t}, \upsilon_{j,t}$] Voltage amplitude of bus $j$ and its square
	\item[$p_{j,t}, q_{j,t}$] Active/reactive power injections at bus $j$ \vspace{1pt}
	\item[$p_{j,t}^{\text{L}}, q_{j,t}^{\text{L}}$] Active/reactive loads of all ELZs at bus $j$ \vspace{1pt}
	\item[$P_{j,t}^{\text{WT}}, Q_{j,t}^{\text{WT}}$] Active/reactive power of WT at bus $j$ \vspace{1pt}
	\item[$P_{j,t}^{\text{PV}}, Q_{j,t}^{\text{PV}}$] Active/reactive power of PV at bus $j$ \vspace{1pt}
	\item[$P_{j,t}^{\text{G}}, Q_{j,t}^{\text{G}}$] Active/reactive power of the grid at bus $j$ \vspace{1pt}
	\item[$k_{j'j,t},\delta_{j'j,k,t}$] OLTC turns ratio and its binary decision variable for the $k$th tap position of NWT $j'j$
\end{IEEEdescription}

\vspace{-8pt}
\subsection{Parameters}
\begin{IEEEdescription}[\IEEEusemathlabelsep\IEEEsetlabelwidth{superscript}]
\addcontentsline{toc}{section}{Nomenclature}

%\item[${N}_s$] Number of representative scenarios
\item[${N_t}, \Delta t$] Scheduling horizon and step length
\item[$\Delta \tau$] Real-time dispatch step length
\item[${N}$] Number of ELZs
\item[$c^{\text{H}_2}$] Hydrogen price
\item[$c^{\text{SU}}, c^{\text{SD}}$] Startup/shutdown costs of ELZ
\item[$c^{\text{G}}$] Electricity purchase price
%\item[$\overline{W}^{\text{C}}$] Upper limit of the capacity of var compensation
\item[$\overline I, \underline I$] Electrolytic current limits of ELZ
\item[$\overline T, \underline T$] Stack temperature limits of ELZ
\item[$K^\text{rec}_\text{min},\Delta K$] Minimum RCT turns ratio and the tap step
\item[$k^{\text{rec}}_\text{max}$] Maximum OLTC tap position
\item[$S^\text{sc}$] PCC short-circuit capacity
\item[$S^\text{GB}, I^\text{GB}$] Base short-circuit capacity and harmonic limit
\item[$\eta^{\text{Cool}}$] Cooling efficiency of ELZ
\item[$P^{\text{By}}$] Standby power consumption of ELZ
\item[$C^{\text{ELZ}}, R^{\text{Diss}}$] Heat capacity and dissipation resistance of ELZ
\item[$T^{\text{Am}}$] Ambient temperature
\item[$c^{\text{Cool}}, T^{\text{Cool}}$] Cooling factor and coolant temperature of ELZ
%\item[$r_{ij}, x_{ij}$] Resistance and reactance of branch $ij$
%\item[$\overline{U}_j, \underline{U}_j$] Voltage magnitude limits at bus $j$
%\item[$\overline{I}_{ij}$] Current capacity limit of branch $ij$
\item[$S_{j}^\text{WT/PV}$] WT/PV installation capacities at bus $j$\vspace{1pt}
\item[$\beta$] Power factor angle limit of PV and the grid
\item[$\overline{k}_{ij}^{\text{All}}$] OLTC action limit of NWT $ij$
\item[$\overline{k}_{ij},\underline{k}_{ij}$] OLTC ratio limits of NWT $ij$
%\item[$U^{\text{rev}},U^{\text{tn}}$] The reversible voltage and thermal neutral voltage
%\item[$N^{\text{cell}},A$] The cell number and electrode area of the stack
%\item[$r_1,r_2,s_1$]  Constant coefficients of U-I characteristics
%\item[$t_1,t_2,t_3 $] Constant coefficients of U-I characteristics
%\item[$F$] The Faraday constant
%\item[$f_1,f_2$] The coefficients of Faradaic efficiency
\end{IEEEdescription}

\section{Introduction}
\label{sec:intro}

\subsection{Background and Motivation}

\IEEEPARstart{R}{enewable} power-to-hydrogen (ReP2H) systems convert variable wind and solar power into low-carbon hydrogen and are increasingly deployed to support decarbonization across energy-intensive sectors \cite{li2025redesigning}. Utility-scale ReP2H plants typically employ multiple electrolyzers (ELZs), whose unit commitment (UC) and load allocation provide flexibility to accommodate renewable generation \cite{varela2021modeling, qiu2023extended}. Because ELZs require high-current DC power, rectifiers form the main electrical interface between the P2H plant and the AC network \cite{ruuskanen2020power, keddar2022power, gao2024advanced}.

Thyristor rectifiers (TRs) are attractive for industrial electrolysis because of their high power rating, efficiency, and low cost \cite{zhang2027multi, gao2024advanced}. For example, the \textit{Sinopec Xinjiang Kuqa Green Hydrogen Demonstration Project} \cite{NARI_Kuqa_2023} and the \textit{China Energy Engineering Songyuan Hydrogen Industry Park} \cite{zeng2025optimal} employ 52 and 32 TRs, respectively. However, their phase-controlled rectification produces characteristic current harmonics at the point of common coupling (PCC), which must satisfy standards such as IEEE 519 \cite{ieeestd} and China's GB/T 14549-1993 \cite{GBT14549}. IGBT-based rectifiers offer better harmonic performance but at substantially higher cost \cite{zeng2025optimal}.

\begin{figure}[t]
	\centering
	\includegraphics[width=3.5in]{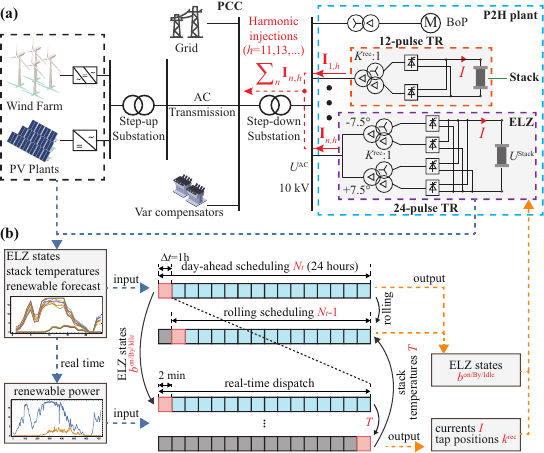}\vspace{-2pt}
	\caption{(a) Structure of a typical grid-connected ReP2H system. (b) Day-ahead/rolling scheduling and real-time dispatch coordination
		framework.}
	\label{fig:system}
\end{figure}

For TR-based systems, harmonics are commonly mitigated through multi-pulse rectification or filters, both of which increase hardware cost and complexity.
12-pulse TRs (12-TRs) use dual-secondary $\Delta$--Y/$\Delta$--$\Delta$ transformer connections to provide a $30^\circ$ phase displacement to mitigate harmonics at relatively low cost, as shown in Fig. \ref{fig:system}(a). Higher-pulse configurations, such as 24-TRs, use additional phase-shifting transformers and bridges to suppress 11th-, 13th-, and higher-order harmonics, but these transformers can account for approximately 40\% of the rectification-stage investment. This motivates a system-level harmonic mitigation (HM) approach that allows low-cost 12-TRs to meet harmonic limits without higher-pulse hardware.

Such an approach is possible because harmonic currents are phasors \cite{zeng2026harmonic}. Harmonics generated by multiple TR-fed ELZs superpose at the PCC, while their magnitudes and phases vary with electrolytic current and low-voltage-side AC voltage. ELZ UC, current allocation, and on-load tap changer (OLTC) control of rectifier transformers (RCTs) can therefore be coordinated to optimize hydrogen production while maintaining harmonic compliance and AC-network voltage security.

\subsection{Literature Review and Research Gap}

ReP2H operation has been extensively studied considering conversion efficiency \cite{wang2024optimization}, ELZ temperature dynamics \cite{ding2024electrolyzer, qiu2023extended, aguado2025optimizing}, pressure \cite{aguado2025optimizing}, hydrogen-in-oxygen (HTO) impurity \cite{qiu2023extended}, and stack degradation \cite{wang2025collaborative, tang2025optimization}. Wang et al. \cite{wang2024optimization} coordinated ELZ commitment and load allocation to improve hydrogen production efficiency. Ding et al. \cite{ding2024electrolyzer} incorporated temperature-dependent operating limits into system optimization, while Aguado et al. \cite{aguado2025optimizing} optimized ELZ temperature and pressure under varying loads. Qiu et al. \cite{qiu2023extended} modeled HTO accumulation and scheduled ELZs to widen their feasible load range. Wang et al. \cite{wang2025collaborative} further incorporated degradation into rolling optimization. However, these studies mainly optimize hydrogen production, with limited attention to grid-side power quality.

For harmonic mitigation, existing methods mainly include multi-pulse rectification \cite{puteanus2024multipulse}, active and passive filters \cite{singh2007improved}, hybrid rectifier topologies \cite{meng2021novel}, and DC-side current shaping \cite{wang2026harmonic}. For example, \cite{meng2021novel} added an auxiliary converter to suppress DC ripple and AC harmonics, while \cite{wang2026harmonic} used a controllable interphase reactor to reshape the DC current of 12-pulse TRs. These device-level methods are effective but require additional hardware or specialized converter designs, increasing the cost and complexity of multi-ELZ plants.

An alternative is to exploit harmonic cancellation among parallel converters. Previous studies coordinated the phases of multiple parallel current sources to reduce aggregated harmonics \cite{yang2026dc, yang2016enhanced}. Our preliminary work \cite{zeng2026harmonic} extended this principle to multi-ELZ plants and demonstrated harmonic cancellation through electrolytic current allocation. The harmonic model and cancellation mechanism were validated using electromagnetic transient simulations in MATLAB/Simulink \cite{zeng2026harmonic}. However, plant-level implementation also requires AC-network voltage security because ELZ loading and RCT tap adjustment jointly affect active and reactive power demand \cite{zeng2024scheduling}. To enable cost-effective harmonic-compliant operation of TR-based ReP2H systems, three practical gaps remain:

\begin{itemize}
	\item RCT tap positions affect firing angles, commutation overlaps, and harmonic phasors, providing an additional degree of freedom for HM. Still, their coordination with electrolytic currents has not been fully modeled or exploited.

	\item Harmonic compliance and voltage security must be maintained continuously at the minute level \cite{GBT14549}, requiring coordinated renewable power tracking and HM.
	
	\item Whether current and tap coordination can make 12-TRs a practical harmonic-compliant alternative to 24-TRs remains unclear.
\end{itemize}

These gaps call for a plant-level method that jointly considers hydrogen production, harmonic compliance, and AC-network voltage security. ELZ UC, load allocation, and RCT tap control must therefore be coordinated across scheduling and real-time dispatch under variable renewable generation.

\subsection{Contributions}

To address these gaps, this paper proposes a two-layer operation scheme for ReP2H systems supplied by 12-TRs, embedding harmonic compliance into plant production scheduling and real-time dispatch. The main contributions are:

\begin{enumerate}
	\item
	A comprehensive harmonic model of TR-fed ELZs is developed to quantify harmonic phasors as functions of AC voltage, RCT tap position, electrolytic current, and stack temperature. The model reveals how coordinated current allocation and RCT tap adjustment across ELZs can satisfy grid-code limits.
	
	\item   A two-layer hydrogen production and HM framework is proposed. Hourly scheduling determines day-ahead and rolling ELZ commitment using the harmonic feasible region (HFR), while minute-level dispatch allocates electrolytic currents for renewable power tracking and adjusts RCT taps for HM.
	
	\item Case studies on real-life systems show that the proposed method (PM) reduces the harmonics of 12-TRs by 48\% and increases profit by 31\% over current-only regulation. Annual simulations further show that coordinated 12-TRs can provide a cost-effective harmonic-compliant alternative to 24-TRs, reducing rectification-stage cost by 37.5\%.
\end{enumerate}

The remainder of this paper is organized as follows. Section \ref{sec:2} develops the ELZ harmonic model. Section \ref{sec:3} presents the two-layer operation scheme. Section \ref{sec:cases} provides case studies, and Section \ref{sec:discussion} discusses engineering implications. Section \ref{sec:conclusion} concludes the paper.

\vspace{2pt}
\section{Harmonic Model of Multi-ELZ P2H Plants}
\label{sec:2}

Fig. \ref{fig:system}(a) shows the grid-connected ReP2H system, comprising wind/photovoltaic (PV) generation, a P2H plant, and an AC network connected to the grid at the PCC. Each ELZ is supplied by a 12-TR through an RCT. The plant consumes active and reactive power while injecting characteristic harmonic currents, whose phasors superpose at the PCC and must satisfy grid-code limits. Without coordination, the current THD of a 12-TR exceeds 6\% \cite{gao2024advanced, Wu2017Multipulse}, compared with approximately 3\%--5\% for a 24-TR \cite{Wu2017Multipulse}. To achieve compliance without expensive higher-pulse hardware, this paper coordinates electrolytic currents and RCT OLTCs for harmonic cancellation, while network OLTCs and var compensators maintain voltage security.

\begin{figure}[t]
	\centering
	\includegraphics[width=3.5in]{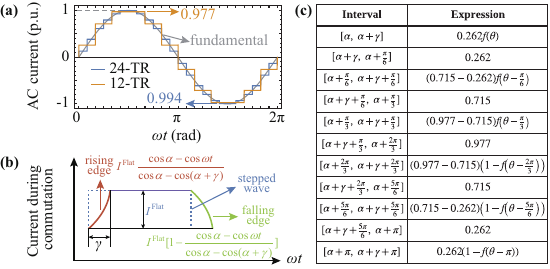}\vspace{-6pt}
	\caption{(a) The rectangular AC current waveforms of the 12- and 24-pulse TRs. (b) The AC current during commutation. (c) Positive half-cycle expression $i^+(\theta)$ of the AC current for a 12-TR.}
	\label{fig:24}
\end{figure}

\subsection{Harmonic Model of Multi-Pulse TR}
\label{sec:harmmodel}

We first relate TR harmonic currents to the firing angle $\alpha$ and commutation overlap $\gamma$. Ideally, the AC-side currents of 24- and 12-pulse TRs are stepped rectangular waveforms, as shown in Fig. \ref{fig:24}. With commutation inductance, the rising $I^\text{Rise}$ and falling $I^\text{Fall}$ intervals are illustrated in Fig. \ref{fig:24}(b) and expressed as
\begin{align}
	&f(\theta)=({\cos\alpha-\cos \theta})/[{\cos\alpha-\cos(\alpha+\gamma)}],\\
	&I^\text{Rise} = I^\text{Flat}f(\theta), ~I^\text{Fall} = I^\text{Flat}\big[1-f(\theta)\big], \label{eq:pH}
\end{align}
where $I^\text{Flat}$ is the current during the flat interval. The positive half-cycle current $i^+(\theta)$ of a 12-TR is shown in Fig. \ref{fig:24}(c), and the 24-TR waveform can be derived similarly. Fourier analysis gives the normalized $h$th harmonic phasor
$\hat{\textbf{I}}_h(\alpha,\gamma) = {\textbf{I}_h}/{I_\text{1st}}$.

\subsection{Harmonic Characteristics of ELZs}

Eqs.  (\ref{eq:Ustack})--(\ref{eq:pstack})  link $\textbf{I}_h(\alpha,\gamma)$ to the electrolytic current $I$, stack temperature $T$, RCT turns ratio $K^\text{rec}$, and AC-side voltage $U^\text{AC}$ through $\alpha$, $\gamma$, and $I_\text{1st}$:
\begin{align}
	&U^{\text{Stack}}(I,T)=N^{\text{Cell}} \big[ U^{\text{rev}}+({r_{1}+r_{2}T})I/{A} \nonumber \\
	& ~~~~~~~~~~~~~~~~~+s_1\log\left({({t_{1}+t_{2}/T+t_{3}/T^{2}})}I/{A} + 1\right) \big],  \label{eq:Ustack}\\
	& U^{\text{Stack}}=U^{\text{AC}}({2.422}/{K^\text{rec}})\cos\alpha-\Delta U, \label{eq:Ustack1}\\
	&\Delta U=U^{\text{AC}}\frac{2.422}{K^\text{rec}}\frac{\cos\alpha-\cos(\alpha+\gamma)}{2}=\frac{3}{\pi} X_\text{c} I, \label{eq:deltaU}\\
	& \sqrt{3}U^{\text{AC}}I_{\text{1st}} \big[{\cos\alpha+\cos(\alpha+\gamma)}\big]/{2}=P^{\text{Stack}}+P^{\text{Loss}}, \label{eq:I1st}\\
	& P^{\text{Stack}}=U^{\text{Stack}}I,~P^{\text{Loss}} = a_2^{\text{Loss}}I^2+a_1^{\text{Loss}}I+a_0^{\text{Loss}}, \label{eq:pstack}
\end{align}
\noindent
where (\ref{eq:Ustack})--(\ref{eq:Ustack1}) describe the stack voltage $U^{\text{Stack}}$, commutation voltage drop $\Delta U$, and $I$ \cite{ulleberg2003modeling}; (\ref{eq:deltaU}) gives $\Delta U$, with 2.422 denoting the rectification coefficient of the 12-pulse TR; ${I}_\text{1st}$ is determined via power conservation (\ref{eq:I1st})--(\ref{eq:pstack}) \cite{li2024two}.
In (\ref{eq:Ustack}), $N^{\text{Cell}}$ is the number of cells in the stack; $U^{\text{rev}}$ is the reversible voltage; $A$ is the electrode area; and $r_1$, $r_2$, $s_1$, $t_1$, $t_2$ and $t_3$ are constants. In (\ref{eq:deltaU}), $X_\text{c}$ is the commutation reactance. In (\ref{eq:pstack}), $a_2^{\text{Loss}}$, $a_1^{\text{Loss}}$ and $a_0^{\text{Loss}}$ are loss-model coefficients.
Solving \eqref{eq:Ustack}--\eqref{eq:deltaU} yields
\begin{align}
	\alpha &= \arccos[(U^{\text{Stack}}+{3}/{\pi} \times X_\text{c}I) {K^\text{rec}}/{(2.422U^{\text{AC}})}], \\
	\gamma &= \arccos[(U^{\text{Stack}}-{3}/{\pi} \times X_\text{c}I) {K^\text{rec}}/{(2.422U^{\text{AC}})}]-\alpha.
\end{align}
% while the coefficient for the 24-TR is 2.4425 \cite{zeng2026harmonic}.

Substituting $\alpha$ and $\gamma$ into $\hat{\textbf{I}}_h(\alpha,\gamma)$ yields the complete expression of $\textbf{I}_h(I,T,U^\text{AC},K^\text{rec})$.
Detailed harmonic characteristics with respect to $I$ are provided in \cite{zeng2026harmonic} and not repeated here.

\begin{figure}[t]
	\centering
	\includegraphics[width=3.5in]{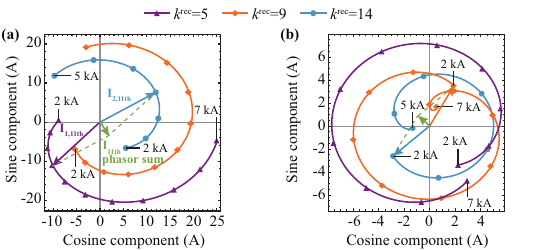}\vspace{-4.5pt}
	\caption{Effect of rectifier OLTC adjustment on (a) 11th and (b) 23rd harmonic currents under different operating currents $I$.}
	\label{fig:cancel}
\end{figure}

In application, each RCT is equipped with a wide-range vacuum OLTC, originally designed to maintain the firing angle within a preferable range under varying load. Here, its tap regulation is exploited for harmonic cancellation among ELZs. The turns ratio $K^\text{rec}$ is a discrete tap-dependent variable:
\begin{align}
	K^\text{rec}=K^\text{rec}_\text{min}+ k^\text{rec} \Delta K,~ 0\leq k^\text{rec}\leq k^\text{rec}_\text{max},~ k^\text{rec}\in Z. \label{eq:krec}
\end{align}
%where the operating range of $k^\text{rec}$ is determined according to the existing TR technical scheme.

\subsection{Harmonic Cancellation via OLTC Adjustment and Harmonic Feasible Region}
\label{sec:harmcha}

Taking two ELZs as an example, Fig. \ref{fig:cancel} illustrates harmonic cancellation through OLTC adjustment at different operating currents $I$. After a tap change, the built-in current controller adjusts $\alpha$ to track the electrolytic-current reference. Changing the RCT tap ratio alters both the magnitude and phase of harmonic currents, creating flexibility for phasor cancellation. For example, at $I=3.5$ kA, setting the two ELZ taps to $k^\text{rec}=14$ and $k^\text{rec}=5$ reduces the 11th harmonic by more than 80\% compared with the centered setting $k^\text{rec}=9$.

%As the studied system is based on a grid-connected ReP2H project in China, GB/T 14549-1993 is adopted as an example standard to formulate the harmonic-compliance constraints at the PCC.
GB/T 14549-1993 \cite{GBT14549} is adopted here to formulate the PCC harmonic-current limits. For an $N$-ELZ P2H plant, the $h$th harmonic must satisfy
\begin{align}
	\big| \sum\nolimits_{n=1}^N {\textbf{I}}_{n,h} \big|  \le \overline{I}_h=(S^{\text{sc}}/S^{\text{GB}})I_h^{\text{GB}}.
\end{align}

Electrolytic currents and OLTC taps are therefore coordinated to modulate the aggregated harmonic phasor. The $N$ ELZs are divided into two-ELZ groups, each assigned a harmonic limit of $2\overline{I}_h/N$, with $N$ typically even in practical projects \cite{zeng2025optimal}. This equal allocation gives a conservative sufficient condition for plant-level compliance while keeping group-wise HM independent and facilitating parallel real-time dispatch.

\begin{figure}[t]
   \includegraphics[width=3.5in]{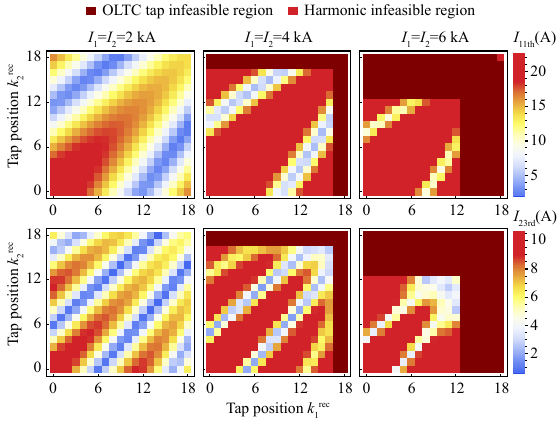}\vspace{-6pt}
  \caption{Harmonic tests under different electrolytic currents and OLTC taps.}
  \label{fig:2K}\vspace{8pt}
  \includegraphics[width=3.35in]{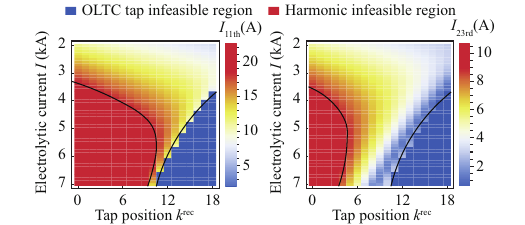}\vspace{-6pt}
  \caption{Harmonic currents of a single ELZ under different operating currents $I$ and RCT tap positions $k^\text{rec}$.}
  \label{fig:1ELZ}
\end{figure}

Fig. \ref{fig:2K} maps the harmonic characteristics under different electrolytic currents and OLTC positions. The HFR contains operating points that satisfy the assigned limits, whereas the OLTC-infeasible region corresponds to invalid firing angles ($\alpha < \alpha_\text{min}=5^\circ$, which reserves sufficient margin for reliable firing) caused by AC/DC voltage mismatch. Even when identical ELZ currents lead to the largest aggregated harmonic injection \cite{zeng2026harmonic}, changing the taps modifies $\alpha$ and $\gamma$ and can reduce the resulting phasor sum.

\begin{remark}
	The harmonic analysis shows that OLTC adjustment can prevent harmonic violations for arbitrary current allocations within a two-ELZ group. Because tap changes alter $\alpha$ and $\gamma$ while maintaining the electrolytic current, they also affect the TR's reactive load. This coupling is included in the AC-network voltage and reactive power models in Sections \ref{sec:re} and \ref{sec:real-time}.
\end{remark}

Fig. \ref{fig:1ELZ} shows the harmonics of a single ELZ under different currents and RCT tap positions. When only one ELZ in a group is online, its harmonic injection remains within the prescribed limit \cite{zeng2026harmonic}. At high load, however, the feasible region for the $(12\pm1)$th harmonics becomes narrow, with only one or two tap positions satisfying the constraints. To improve robustness against real-time short-circuit capacity (SCC) variations, a safety margin is imposed on the single-ELZ HFR to avoid violations under weaker grid conditions:
\begin{gather}
	b^{\text{On}}\underline I\leq I\leq b^{\text{On}}\overline I,  \label{eq:I}\\
	I_1 \leq \overline{I}^{\text{single}}+ (\overline{I}-\overline{I}^{\text{single}})b_2^{\text{On}},\\
	I_2 \leq \overline{I}^{\text{single}}+(\overline{I}-\overline{I}^{\text{single}})b_1^{\text{On}}, \label{eq:I3}
\end{gather}
where $\overline{I}^{\text{single}}$ is the current limit when the other ELZ in the group is idle. In this study, $\overline{I}^{\text{single}}$ is set to 6 kA.

\section{Two-Layer Coordination of Production Scheduling and Harmonic Mitigation}
\label{sec:3}

Beyond harmonic compliance, operation of TR-based ReP2H systems must consider not only hydrogen production but also the TRs' reactive load and its impact on voltage and network losses. This section therefore embeds HM into plant operation through a two-layer framework, in which hourly scheduling determines ELZ states within the HFR under renewable uncertainty, while real-time dispatch tracks renewable power and coordinates HM.

\subsection{Operational Model of the ReP2H System}

\subsubsection{Network Power Flow}

The ReP2H system generally has a radial topology \cite{zeng2024scheduling}; therefore, the DistFlow model \cite{zeng2024scheduling, farivar2013branch} is used to represent branch power flows, losses, and voltage profiles, following (26)--(32) in \cite{zeng2024scheduling}. Nodal injections include wind/PV generation, grid exchange, var compensation, and the $N$-ELZ P2H load. The network thus follows:
\begin{align}
	& \text{DistFlow branch power flow (26)--(32) in \cite{zeng2024scheduling}}, \label{eq:29}\\
	&p_{j,t}=P_{j,t}^{\text{WT}}+P_{j,t}^{\text{PV}}+P_{j,t}^{\text{G}}-p_{j,t}^{\text{L}} ,\label{eq:36}\\
	&q_{j,t}=Q_{j,t}^\text{WT}+Q_{j,t}^\text{PV}+Q_{j,t}^\text{C}+Q_{j,t}^\text{G}-q_{j,t}^{\text{L}},  \label{eq:37}\\
	&p_{j,t}^{\text{L}}=\textstyle\sum_{n=1}^{{N}} P_{n,t}^{\text{ELZ}}, \ q_{j,t}^\text{L}=\textstyle\sum_{n=1}^{{N}}Q_{n,t}^\text{ELZ} . \label{eq:39}
\end{align}

\subsubsection{WTs and PV Plants}

The active and reactive powers of WTs and PV plants satisfy:
\begin{align}
	& 1.24P_{j,t}^\text{WT} - 0.91S_{j}^\text{WT}\le Q_{j,t}^\text{WT} \le0.91S_{j}^\text{WT}-0.58P_{j,t}^\text{WT}, \label{eq:qwt}\\
	& \left(Q_{j,t}^{\text{PV}}\right)^2\le\left(S_j^{\text{PV}}\right)^2-\left(P_{j,t}^{\text{PV}}\right)^2,~ |Q_{j,t}^{\text{PV}}|\le P_{j,t}^{\text{PV}}\tan\beta, \label{eq:qpv}
\end{align}

\subsubsection{Interaction with the Grid}

The P2H plant is allowed to purchase electricity from the grid, but the power factor at the PCC must not violate the prescribed limit \cite{li2024two}:
\begin{align}
	P_{j,t}^{\text{G}}\geq0, ~|Q_{j,t}^{\text{G}}|\le P_{j,t}^{\text{G}}\tan\beta,
\end{align}

\subsubsection{Var Compensation and OLTCs}

The reactive power compensation provided by static var generators (SVGs) and OLTC at the step-down substation is modeled as follows \cite{ding2017data}:
\begin{align}
	&\begin{vmatrix}Q_{j,t}^\text{C}\end{vmatrix}\le\overline{Q}_{j}^\text{C}, \hspace{110pt}\label{Qc}
\end{align}
\vspace{-16pt}

\break
\
\vspace{-12pt}
\begin{align}
	&\begin{cases}
		k_{j'j,t}=\sum_{k=0}^{K_{j'j}}\delta_{j'j,k,t}w_{j'j,k},\\
		\sum_{k=0}^{K_{j'j}}\delta_{j'j,k,t}=1,\delta_{j'j,k,t}\in\{0,1\}
	\end{cases}  \label{eq:k=1}\\
	&\begin{cases}
		\sum_{t=2}^{{N_t}}\left|k_{j'j,t}-k_{j'j,t-1}\right|\le\overline{k}_{j'j}^{\text{All}},\\
		\underline{k}_{j'j}\le k_{j'j,t}\le\overline{k}_{j'j},
	\end{cases} \label{eq:k}
\end{align}

\subsubsection{State Switching of ELZs}

Each ELZ operates in one of three mutually exclusive states, i.e., \textit{production}, \textit{standby}, or \textit{idle}. Startup, shutdown, and transition delays are represented by \cite{varela2021modeling}:
\begin{align}
	&b_{n,t}^{\text{On}}+b_{n,t}^{\text{By}}+b_{n,t}^{\text{Idle}}=1, \label{eq:logic}\\
	&b_{n,t}^{\text{On}}+b_{n,t}^{\text{By}}+b_{n,t-1}^{\text{Idle}}-1\leq b_{n,t}^{\text{SU}}, \label{eq:1}\\
	&b_{n,t-1}^{\text{On}}+b_{n,t-1}^{\text{By}}+b_{n,t}^{\text{Idle}}-1\le b_{n,t}^{\text{SD}}, \label{eq:2}\\
	-&b_{n,t-2}^{\text{Idle}}-b_{n,t}^{\text{Idle}}+b_{n,t-1}^{\text{Idle}} \leq 0. \label{eq:3}
\end{align}

\subsubsection{Active Power and Hydrogen Yield of ELZs}

The active power of each ELZ consists of electrolytic power, balance-of-plant (BoP) consumption, and rectifier losses \cite{li2024two}. Hydrogen  yield depends on the stack current and Faradaic efficiency \cite{ulleberg2003modeling}.
\begin{align}
	&P_{n,t}^\text{ELZ}=P_{n,t}^\text{Stack}+P_{n,t}^\text{BoP}+P_{n,t}^\text{Loss}, \label{eq:Pele} \\
	&P_{n,t}^{\text{BoP}}=\left(1-b_{n,t}^{\text{Idle}}\right)P_{n,t}^{\text{Cool}}/\eta^{\text{Cool}}+b_{n,t}^{\text{By}}P^{\text{By}}, \label{eq:Paux}\\
	%   &P_{n,t}^{\text{Loss}} = a_2^{\text{Loss}}I_{n,t}^2+a_1^{\text{Loss}}I_{n,t}+a_0^{\text{Loss}}, \label{eq:Ploss}\\
	&Y_{n,t}^{\text{H}_2}=3.6 \eta_{n,t}^{\text{F}}N^{\text{Cell}}I_{n,t}/F, \label{eq:YH2}\\
	&\eta_{n,t}^{\text{F}}={(I_{n,t}/A)^{2}}/[{f_{1}+(I_{n,t}/A)^{2}}] \times f_{2}, \label{eq:etaF1}
\end{align}
where $F$ is the Faraday constant; $\eta_{n,t}^{\text{F}}$ is the Faradaic efficiency; $f_{1}$ and $f_{2}$ are the coefficients of Faradaic efficiency.

%\subsubsection{Degradation of Electrolytic Stack}
%The degradation of industrial alkaline ELZs is mainly reflected in increased stack voltage $\Delta U^\text{Stack}$ and energy consumption, driven by thermal/mechanical stress from start-stop cycles, and local overpotentials and thermal stress caused by power fluctuations, and is modeled as follows \cite{wang2025collaborative, zheng2023off}:
%\begin{align}
%  \Delta U^\text{Stack}=\Delta U^\text{SU} b^\text{SU}+\Delta U^\text{SD} b^\text{SD}+\Delta U^\text{adj} \dfrac{|\Delta P^\text{Stack}|}{P^\text{Stack}_\text{max}}, \label{eq:degradation}
%\end{align}
%where $\Delta U^\text{SU}$, $\Delta U^\text{SD}$, and $\Delta U^\text{adj}$ are the relevant coefficients.

\subsubsection{Reactive Power Demand of ELZs}
\label{sec:re}

Because HM changes the firing and overlap angles, it also affects TRs' reactive power and hence network voltage and losses. The reactive power of each ELZ includes fundamental phase-shift and harmonic-distortion components \cite{zeng2024scheduling,keddar2022power}:
\begin{align}
	&Q^{\text{ELZ}} = %\big\Vert{({Q^{\text{S}}}, {Q^{\text{D}}})} \big\Vert _2,
	\sqrt{({Q^{\text{S}}})^{2}+({Q^{\text{D}}})^{2}},
	\label{eq:Qstack} \\
	&Q^{\text{S}}=\sqrt{3}U^{\text{AC}}I_{\text{1st}}\sin\varphi,~ Q^{\text{D}}=\sqrt{3}U^{\text{AC}}I_{\text{1st}}{\sqrt{1-\nu^{2}}}/{\nu}. \label{eq:Qd}\\
	&\varphi=\arccos \big[({\cos\alpha+\cos(\alpha+\gamma)})/{2}\big],
\end{align}
where $\nu$ is the harmonic factor; and $\varphi$ is the power factor angle.

\subsubsection{Thermal Dynamics of ELZs}
ELZ temperature is governed by heat generation, dissipation, and cooling, as:
 \begin{align}
  &C_n^{\text{ELZ}}\left(T_{n,t+1}-T_{n,t}\right)=\left(P_{n,t}^{\text{Heat}}-P_{n,t}^{\text{Diss}}-P_{n,t}^{\text{Cool}}\right)\Delta t, \label{eq:thermal}\\
  &P_{n,t}^{\text{Heat}}=P_{n,t}^{\text{Stack}}-N^{\text{Cell}}U^{\text{tn}}I_{n,t},  \label{eq:Pgen}\\
  &P_{n,t}^{\text{Diss}}=({T_{n,t}-T^{\text{Am}}})/{R_n^{\text{Diss}}},\\
  &0\leq P_{n,t}^{\text{Cool}} \leq c^{\text{Cool}}\left(T_{n,t}-T^{\text{Cool}}\right),~ \underline T\leq T \leq \overline T \label{eq:T}
 % &\underline T\leq T \leq \overline T.  \label{eq:T}
\end{align}
where $U^{\text{tn}}$ is the thermal neutral voltage.

More details of the operational model are provided in our previous work \cite{zeng2024scheduling}. The nonlinear terms are handled by polynomial fitting, piecewise linearization, and Big-M reformulations, so that the overall model is transformed into a mixed-integer second-order cone programming (MISOCP) problem.

\subsection{Rolling Scheduling and Real-Time Dispatch Coordination}

Fig. \ref{fig:system}(b) shows the proposed two-layer  framework to coordinate ELZ operation, renewable power tracking, and HM. The upper layer performs hourly rolling scheduling over the remaining horizon. At each update, only the ELZ commitment for the first hour is implemented and passed to the lower-level real-time layer; subsequent decisions are reoptimized in the next scheduling step.

Following GB/T 14549-1993 \cite{GBT14549}, a 2-min resolution is adopted for real-time dispatch and harmonic assessment. Given actual renewable generation and ELZ temperatures, the lower layer allocates electrolytic currents to the ELZs for optimizing power tracking and hydrogen production while adjusting RCT tap positions for HM. Temperatures are updated after each interval. At the end of each hour, the terminal states are returned to the upper layer for the next rolling optimization.

\subsection{Day-Ahead/Rolling Scheduling}
\label{sec:scheduling}

\subsubsection{Two-Stage Stochastic Programming}

The upper-layer scheduling model minimizes the net operating cost $C$, comprising startup/shutdown and electricity-purchase costs minus hydrogen sales revenue. First, a deterministic scheduling model is given as follows:
\begin{align}
	\nonumber \min~& C=\sum_{t=1}^{N_t}\Big[\sum_{n=1}^{N}(c^{\text{SU}}b_{n,t}^{\text{SU}}+c^{\text{SD}}b_{n,t}^{\text{SD}}-c^{\text{H}_2}Y_{n,t}^{\text{H}_2})+c^{\text{G}}P_{t}^{\text{G}}\Big] \Delta t ,\\
	\nonumber \text{s.t. }~& \text{HFR }(\ref{eq:I})-(\ref{eq:I3}),~\text{AC-network operation }(\ref{eq:29})-(\ref{eq:k}), \\
	&\text{ELZ operation }(\ref{eq:logic})-(\ref{eq:T}). \label{eq:all}
\end{align}
Here, maximizing the revenue also naturally promotes renewable energy utilization and reduces network losses. In hourly scheduling, RCT taps are fixed at their centered positions to estimate P2H reactive power; their actual settings and reactive-power demand are updated in real-time dispatch.

To address renewable uncertainty, a stochastic programming model is further developed. Renewable uncertainty is represented by scenarios generated from forecast-error distributions using Latin hypercube sampling and reduced to $N_s$ representative scenarios with probabilities $\pi_s$ \cite{li2024restoration}. Because ELZ commitment must be decided before renewable generation is realized, the problem is formulated in two stages. The first-stage variables $\bm{x}$ represent non-anticipative ELZ commitments, while $\bm{y}_s$ contains scenario-dependent operating decisions:
\begin{align}
	\min_{\{\bm{x},\bm{y}_s,\forall s\}}&\bm{a}^\text{T}\bm{x}+\sum\nolimits_{s=1}^{N_s}\pi_s \bm{b}_s^\text{T}\bm{y}_s,  \label{eq:PHobj}\\
	\text{s.t.}~~~&\bm{x}\in X, \bm{y}_s\in Y_s, \forall s,	  \label{eq:PHst}
\end{align}
where $\bm{a}$ and $\bm{b}_s$ are coefficient vectors; $X$ and $Y_s$ are the feasible sets of $\bm{x}$ and $\bm{y}_s$; $N_s$ is the number of representative scenarios.

\subsubsection{Solution Method}

Due to the large number of scenarios and binary variables, the resulting mixed-integer stochastic problem is difficult to solve.
To address the issue, we use progressive hedging (PH) to decompose the stochastic problem into $N_s$ scenario subproblems (SPs) while enforcing non-anticipativity of ELZ commitment through multiplier and penalty updates \cite{rockafellar1991scenarios}. Numerical tests show that conventional PH may fail to converge when only a few binary commitments remain inconsistent. We therefore adopt a modified PH (MPH) method  \cite{liu2025application}, and the overall solution procedure is presented in Algorithm \ref{alg:MPH}.

PH iterations continue until fewer than three binary variables remain inconsistent, i.e., $z\leq2$. Their possible combinations are then enumerated, and the expected objective is recalculated with these binaries fixed. The best combination is selected. Because at most $2^z\leq4$ cases are evaluated, this final enumeration resolves residual inconsistency with little computational burden.

\begin{algorithm2e}[tb]
	\caption{MPH Algorithm for the Scheduling Model}%\footnotesize
	\label{alg:MPH}
	\DontPrintSemicolon
	\LinesNumbered
	
	Initialize the iteration index $m\gets0$, Lagrange multiplier
	$\bm{\omega}_s^0\gets0$, and index $z^0$\;
	
	Get the continuous and integer variables in parallel by
	$\{\bm{x}_s^0,\bm{y}_s^0\}
	\leftarrow
	\arg\min_{\{\bm{x},\bm{y}_s\}}
	\left\{
	\bm{a}^{\text{T}}\bm{x}
	+\bm{b}_s^{\text{T}}\bm{y}_s:
	\text{(\ref{eq:PHst})}
	\right\},\quad \forall s$
	
	Calculate the initial value
	$\bm{\overline{x}}^0\leftarrow\sum_s\pi_s\bm{x}_s^0$\;

	\While{$z\geq \epsilon~(\epsilon=3)$}{
		
		Update the index $m\gets m+1$\;
		
		\For{$s=1$ \KwTo $N_s$}{
			Update the multiplier
			$\bm{\omega}_s^m
			\gets\bm{\omega}_s^{m-1}+\rho(\bm{x}_s^{m-1}-\bm{\overline{x}}^{m-1})$\;
			
			Update the first- and second-stage variables $\{\bm{x}_s^m,\bm{y}_s^m\}\gets
			\arg\min_{\{\bm{x},\bm{y}_s\}}
			\{\bm{a}^{\text{T}}\bm{x}
			+\bm{b}_s^{\text{T}}\bm{y}_s+(\bm{\omega}_s^m)^{\text{T}}\bm{x}+
			\frac{\rho}{2}\|\bm{x}-\bm{\overline{x}}^{m-1}\|_2^2:
			\text{(\ref{eq:PHst})}
			\}$\;
		}
		
		$\bm{\overline{x}}^m
		\leftarrow
		\sum_s\pi_s\bm{x}_s^m$\;
		
		$z=$ number of non-binary elements in $\bm{\overline{x}}^m$\;
	}
	
	Locate all integers $\bm{\overline{x}}_{\text{ic}}$ violating non-anticipativity constraints, and denote the others as $\bm{\overline{x}}_{\text{c}}^m$\;
	
	Create all distinct combinations of
	$\{\bm{\overline{x}}_{\text{ic}}^1,\ldots,
	\bm{\overline{x}}_{\text{ic}}^{2^z}\}$\;
	
	\For{$case=1$ \KwTo $2^z$}
	{
		$obj^{\text{case}}=\min_{\bm{\overline{x}}_{\text{ic}}^{\text{case}}}\
		\{\bm{a}^{\text{T}}
		\begin{bmatrix}
			\bm{\overline{x}}_{\text{ic}}^{\text{case}}~
			\bm{\overline{x}}_{\text{c}}^m
		\end{bmatrix}
		+\sum_s\pi_s\bm{b}_s^{\text{T}}\bm{y}_s^m:
		\text{(\ref{eq:PHst})}\}$;
	}
	Select the $\bm{\overline{x}}_{\text{ic}}^{\text{case}}$ with the minimum $obj^{\text{case}}$\;
	\Return{$\{\bm{\overline{x}}_{\text{ic}}^\text{case},\bm{\overline{x}}^m,\bm{y}_1^m,\ldots,\bm{y}_{N_s}^m\}$}\;
\end{algorithm2e}

\subsection{Real-Time Dispatch at a 2-min Resolution}
\label{sec:real-time}

The real-time layer needs to jointly optimize hydrogen production (HP) and HM. Directly embedding the highly nonlinear harmonic model into the full network problem would increase computational burden. The problem is therefore decomposed into one HP SP for power tracking and current allocation and $N/2$ independent HM SPs for RCT tap coordination.

\subsubsection{Power Tracking and Allocation}

Given ELZ states and actual renewable output, the HP subproblem allocates electrolytic currents to maximize hydrogen yield subject to network constraints:
\begin{align}
	\min \ &\Big[\sum_{n=1}\nolimits^N(-c^{\text{H}_2}Y_{n,t}^{\text{H}_2})+c^{\text{G}}P_{t}^{\text{G}}\Big] \Delta \tau, \label{eq:50}\\
	\text{s.t.}\ &(\ref{eq:29})-(\ref{eq:k}),~(\ref{eq:Pele})-(\ref{eq:T}). \label{eq:51}
\end{align}

Unlike (\ref{eq:all}), this problem excludes the HFR and ELZ state-switching constraints because commitment is fixed by the upper layer. The resulting electrolytic currents $I$, PCC voltage $U^\text{AC}$, and stack temperatures $T$ are passed to the HM subproblems.

\subsubsection{Harmonic Mitigation}

Each HM subproblem adjusts the electrolytic currents and RCT tap positions to satisfy harmonic limits while minimizing deviation from the production optimum and unnecessary tap changes:
\begin{align}
	&\min_{\{I_1,I_2,k_1^\text{rec},k_2^\text{rec}\}} \ c_1\sum\nolimits_{n=1}^2|I_n-I_n^\text{ref}|+c_2\sum\nolimits_{n=1}^2|\Delta k_n^\text{rec}| \label{eq:52}\\
	&\text{s.t.} ~~ |\sum\nolimits_{n=1}\nolimits^2{\textbf{I}}_{n,h}| \leq 2\overline{I}_h/N,~\alpha_\text{min}\leq\alpha\leq\alpha_\text{max},~(\ref{eq:krec}), \label{eq:Imax}
\end{align}
where $c_1$ is set large to prioritize HP; $c_2$ denotes the RCT OLTC action cost; $\alpha_\text{max}$ is imposed to reserve sufficient commutation margin.
Although (\ref{eq:52})--(\ref{eq:Imax}) form a mixed-integer nonlinear programming (MINLP) problem, each HM subproblem contains only two ELZs and can be solved efficiently.

\subsubsection{Iterative Solution Approach}

Although the HP and HM SPs can be solved sequentially, they remain coupled through reactive power. Changing RCT taps modifies $\alpha$ and $\gamma$, which changes TR reactive power and network voltages; these changes, in turn, affect $I$, $T$, $U^\text{AC}$. Algorithm \ref{alg:ITR} therefore alternates between HM and network-aware HP dispatch until convergence. A limiter is imposed on the update of $U^\text{AC}$ in Line \ref{line:10} to suppress oscillatory iterations.

\begin{algorithm2e}[tb]
	\caption{Iterative HP-HM Solution Approach}%\footnotesize
	\label{alg:ITR}
	\DontPrintSemicolon
	\LinesNumbered
	Initialize $m\gets0$, $k^\text{rec}$, convergence indicator $g^0_\text{I/T/U}$ and convergence tolerance $\epsilon_\text{I/T/U}$\;
	Get the initial $I^0$, $T^0$, and $U^\text{AC,0}$ by solving (\ref{eq:50})--(\ref{eq:51})\;
	\While{$g^m_\text{I/T/U}\ge \epsilon_\text{I/T/U}$}
	{
		Update $m\gets m+1$
		
		\For{$ELZ~group=1$ \KwTo $N/2$}
		{
			Mitigate harmonics by coordinating the electrolytic currents and OLTC taps $\{I_1^m,I_2^m,k_{1}^{\text{rec},m},k_2^{\text{rec},m}\}\gets\arg\min_{\{I_1,I_2,k_1^\text{rec},k_2^\text{rec}\}}
			\{\text{(\ref{eq:52})}:\text{(\ref{eq:Imax})}\}$\;
		}
		Update $K^{\text{rec},m}\gets K^\text{rec}_\text{min}+ k^{\text{rec},m} \Delta K$\;
		Calculate $Q^\text{ELZ}(I^{m-1},T^{m-1},U^{\text{AC},m-1},K^{\text{rec},m})$ and substitute it into (\ref{eq:37})\;
		Track renewable power and optimize HP by solving $\min\{\text{(\ref{eq:50})}:\text{(\ref{eq:51})}, |U^{\text{AC},m}-U^{\text{AC},m-1}|\leq0.01\}$ to obtain $I^{m}$, $T^m$, and $U^{\text{AC},m}$
				\label{line:10}
		\;
		Assess $g^m_\text{I/T/U}\gets |({I/T/U^\text{AC}})^m-({I/T/U^\text{AC}})^{m-1}|$\;
	}
	\Return {$\{I^m, T^m, U^{\text{AC},m}, k^{\text{rec},m}\}$}
\end{algorithm2e}

\begin{figure}[t]
	\centering
	\includegraphics[width=3.48in]{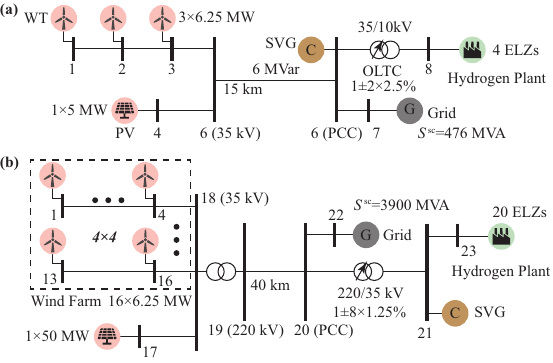}\vspace{-4.5pt}
	\caption{Topologies of the ReP2H systems used for the case studies. (a) Small-scale system with 4 ELZs. (b) Large-scale system with 20 ELZs.}
	\label{fig:case}
\end{figure}

\section{Case Studies}
\label{sec:cases}

\begin{figure}[t]
	\includegraphics[width=3.5in]{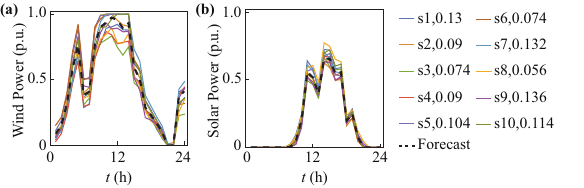}\vspace{-8pt}
	\caption{(a) Wind and (b) PV generation scenarios and their probabilities.}
	\label{fig:WT}
\end{figure}

Case studies are conducted on two ReP2H systems in Inner Mongolia, China, following the configurations in \cite{zeng2025optimal,zhu2026exploring}, as shown in Fig. \ref{fig:case}. The wind and PV scenarios and their probabilities are given in Fig. \ref{fig:WT}. Key parameters are listed in Table \ref{tab:para}, with additional details available in \cite{zeng2025optimal,zeng2024scheduling,zeng2026harmonic}.

\textit{1) Small-scale system:} It comprises $3\times6.25$~MW WTs, a $1\times5$~MW PV plant, a 6~MVar SVG, and four alkaline ELZs. Renewable power is delivered to the P2H plant through a 15-km 35-kV line, and the PCC SCC is 476 MVA. This system is used to demonstrate coordinated HP and HM and assess whether 12-TRs can replace 24-TRs for harmonic compliance.

\textit{2) Large-scale system:} It comprises $16\times6.25$~MW WTs, a $1\times50$~MW PV plant, an SVG, and twenty alkaline ELZs. It is connected to the PCC at 220~kV with an SCC of 3900~MVA. This system evaluates the scalability and computational performance of the PM.

All simulations are implemented in \textit{Wolfram Mathematica 14.0} on a computer with an \textit{Intel Core Ultra 7 165H@1.40 GHz} CPU and 32 GB of RAM.
The MISOCP problems for day-ahead/rolling scheduling and real-time power tracking are solved using \textit{Gurobi 12.0.2}, while the HM SPs are solved with \textit{Mathematica}'s built-in global optimizer. ELZs are paired for HM as (1,2), (3,4), and so forth.

\begin{table}[tb]\scriptsize
	\renewcommand{\arraystretch}{1.2}
	\caption{Key Operational Parameters Used in Case Studies}\vspace{-6pt}
	\label{tab:para}
	\centering
	\begin{tabular}{c@{\hspace{4pt}}c@{\hspace{4pt}}c@{\hspace{4pt}}c}
		\hline\hline
		Parameter                             & Value                                 & Parameter                        & Value   \\
		\hline
		$r_1$                                 & 1.71$\times 10^{-4} \Omega$m$^2$                 & $t_1$                                  & $-$0.24 m$^2$/A           \\
		$r_2$                                 & $-$1.97$\times10^{-7} \Omega$m$^2$/\textcelsius  & $t_2$                                  & 26.23 m$^2$\textcelsius/A            \\
		$s_1$                                 & 0.1618 V                                         & $t_3$                                  & 139.87 m$^2$\textcelsius$^2$/A      \\
		$N^{\text{Cell}}$/$A$                 & 350 / 4.0 m$^2$                                  & $X_\text{c}$                         &0.0072 $\Omega$
		\\
		$\overline I$ / $\underline I$         & 7 / 2 kA      &$\overline T$ / $\underline T$        & 80 / 25 \textcelsius    \\                         $c^{\text{G}}$                      & 0.6  CNY/kWh       &$c^{\text{H}_2}$                      & 22 CNY/kg             \\ $c^{\text{SU}}$ / $c^{\text{SD}}$      & 1000 / 0 CNY    &$c_2$      & 0.5 CNY/action        \\
		$S^\text{GB}$ (35kV)      & 250 MVA    &$I_\text{11/13/23/25}^\text{GB}$ (35kV)     & 5.6 / 4.7 / 2.7 / 2.5 A      \\
		$S^\text{GB}$ (220kV)      &  2000 MVA    &$I_\text{11/13/23/25}^\text{GB}$ (220kV)     & 4.3 / 3.7 / 2.1 / 1.9 A      \\
		$\Delta K$     & $2.5\%$                & $k^\text{rec}_\text{max}$              & 18 \\
		\hline\hline
	\end{tabular}
\end{table}

\subsection{Operational Results of the 4-ELZ Small-Scale System}

\subsubsection{Day-Ahead Scheduling}

\begin{figure}[t]
	\centering
	\includegraphics[width=3.45in]{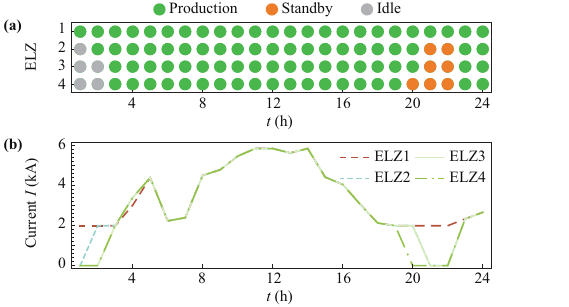}\vspace{-8pt}
	\caption{Day-ahead (a) on-off switching of ELZs and (b) current allocation in one of the ten scenarios.}
	\label{fig:dayahead}
\end{figure}

Fig. \ref{fig:dayahead} shows the day-ahead ELZ on-off switching states and current allocation for one of the ten scenarios. During low renewable generation, some ELZs reduce their currents, remain idle at $t=\{1,2\}$ h, or enter standby at $t=\{21,22\}$ h to avoid grid electricity purchase. During high renewable generation periods, all ELZs are committed and share the load nearly evenly, improving renewable energy utilization and P2H energy conversion efficiency. When only one ELZ in a two-ELZ group is in production, such as around $t=\{1, 20, 22\}$ h, its current remains within the single-ELZ HFR.

\subsubsection{Real-time Dispatch}
\label{sec:real-time-result}

\begin{figure}[t]
	\centering
	\includegraphics[width=3.45in]{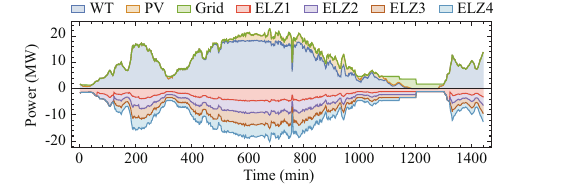}\vspace{-6pt}
	\caption{Power supply and electrolytic load in the real-time layer.}
	\label{fig:power}
\end{figure}

\begin{figure}[t]
	\centering
	\includegraphics[width=3.45in]{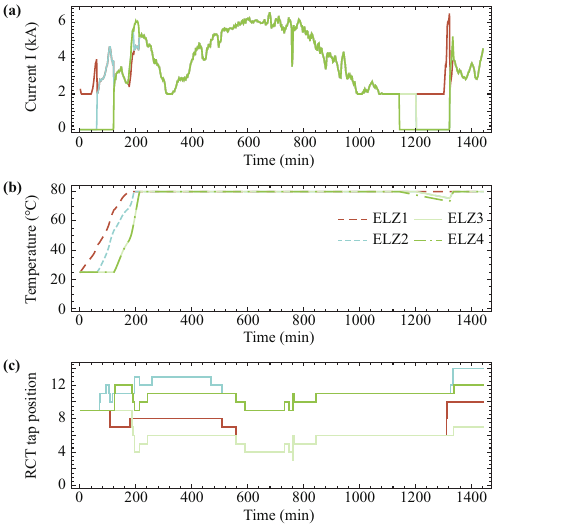}\vspace{-6pt}
	\caption{(a) Electrolytic current, (b) stack temperature, and (c) RCT tap position in real-time dispatch in the 4-ELZ system.}
	\label{fig:IT}
\end{figure}

\begin{figure}[t]
	\centering
	\includegraphics[width=3.45in]{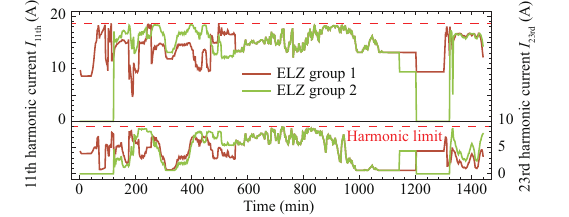}\vspace{-6pt}
	\caption{10-kV-side 11th/23rd harmonic current in the 4-ELZ system.}
	\label{fig:harmonicI}
\vspace{8pt}
	\includegraphics[width=3.45in]{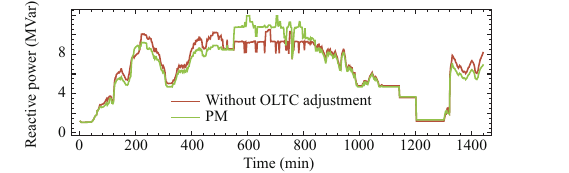}\vspace{-6pt}
	\caption{Reactive power demand of the P2H plant in the 4-ELZ system.}
	\label{fig:repower}
\end{figure}

Fig. \ref{fig:power} shows the real-time balance among renewable generation, grid purchases, and ELZ loads. Fig. \ref{fig:IT} presents electrolytic currents, stack temperatures, and RCT tap positions. The 10-kV-side 11th and 23rd harmonic currents of each ELZ group are shown in Fig. \ref{fig:harmonicI}; similar results for the 13th and 25th harmonics are omitted. Fig. \ref{fig:repower} further presents the reactive power demand of the P2H plant.

As shown in Fig. \ref{fig:power}, the electrolytic loads follow wind and PV variations. Grid power is purchased when renewable generation is insufficient during $\tau \in$ [1080, 1300] min to maintain minimum ELZ loading. Fig. \ref{fig:IT}(a) shows nearly balanced current allocation among online ELZs under the equimarginal principle \cite{zeng2024scheduling}, indicating that HM requires little deviation from the production optimum \cite{zeng2026harmonic}. Instead, most harmonic adjustment is provided by the RCT tap coordination.
Fig. \ref{fig:IT}(c) shows that the two RCT taps in each group are coordinately adjusted to create harmonic-phasor differences, echoing the approximately linear HFR in Fig. \ref{fig:2K}. During high-load intervals, such as $\tau \in$ [520, 800] min, lower RCT tap positions are more favorable for HM, as the feasible tap combinations shift toward lower positions when the electrolytic load increases.

Fig. \ref{fig:harmonicI} confirms that the 11th and 23rd harmonic currents remain below their limits throughout the dispatch horizon. Fig. \ref{fig:repower} shows that the PM can alter the plant reactive power demand compared with operation without OLTC adjustment, while the resulting increase in network losses is only about 0.03\%. All bus voltages remain within 0.954--1.050 p.u., ensuring voltage security. When only one ELZ in a group is online,  its harmonic injection already satisfies the limits,  and unnecessary tap changes are avoided.
These results verify the proposed harmonic-compliant operation without compromising hydrogen production or voltage security.

\subsection{Comparisons With Other Operational Schemes}

To verify the value of coordinating electrolytic current allocation and RCT tap adjustment for harmonic-compliant operation, the following methods are compared:
\begin{itemize}
	\item \textbf{PM}: The proposed operation method coordinating HP and HM, as described in Sections \ref{sec:scheduling} and \ref{sec:real-time}.
	\item \textbf{M1}: Similar to PM, but harmonic constraints are enforced only through electrolytic-current regulation, with fixed RCT taps, following \cite{zeng2026harmonic}.
	\item \textbf{M2}: Harmonic constraints are neglected, and operation optimizes only HP, as commonly assumed in existing studies \cite{ding2024electrolyzer, qiu2023extended,wang2024optimization,wang2025collaborative, tang2025optimization, aguado2025optimizing}.
\end{itemize}

\begin{figure}[t]
	\centering
	\includegraphics[width=3.45in]{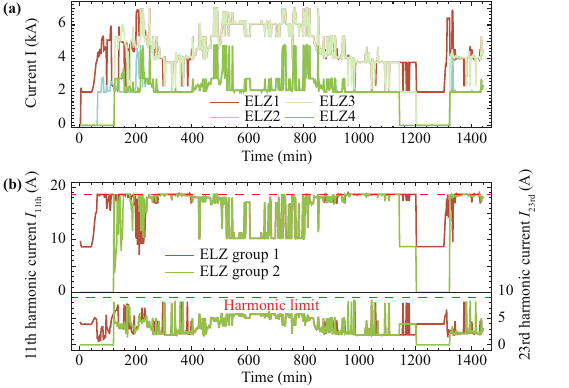}\vspace{-6pt}
	\caption{(a) Electrolytic current and (b) 11th/23rd harmonic current of ELZ groups in the 4-ELZ system under M1.}
	\label{fig:M1}
	\vspace{8pt}
	\includegraphics[width=3.45in]{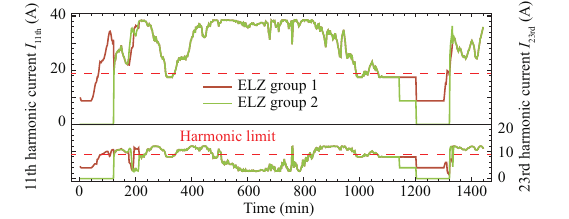}\vspace{-6pt}
	\caption{11th/23rd harmonic current of ELZ groups under M2.}
	\label{fig:M2}
\end{figure}

Table \ref{tab:comparison} compares the performance of PM, M1, and M2, while Figs. \ref{fig:M1} and \ref{fig:M2} show their electrolytic current allocations and harmonic current injections, respectively.

Under M1, HM relies on current differences between the two ELZs in each group, as shown in Fig. \ref{fig:M1}(a) \cite{zeng2026harmonic}. ELZs belonging to different groups, such as ELZ1 and ELZ3, retain similar electrolytic loads to optimize P2H energy conversion efficiency. However, the narrower HFR under current-only regulation limits M1 from fully following power fluctuations, reducing renewable accommodation and lowering hydrogen output.
With comparable and compliant harmonic levels, Table \ref{tab:comparison} shows that PM increases operating profit and hydrogen production by 31.0\% and 13.9\%, respectively, relative to M1.

Without HM, M2 achieves nearly the same hydrogen output as PM but violates harmonic limits. In Fig. \ref{fig:M2}, high-load 11th harmonic exceeds its limit by $\frac{38.5-18.66}{18.66}=$ 106.3\%. Compared to M2, PM causes only a negligible profit reduction due to RCT tap actions, while reducing the average 11th and 23rd harmonic currents by 48.0\% and 43.2\%, respectively. Current and tap coordination thus preserves the production performance of unconstrained operation while restoring harmonic compliance.

\begin{table}[t]
	\renewcommand{\arraystretch}{1.2}
	\caption{Performance Comparison in the 4-ELZ Small-Scale System}\vspace{-6pt}
	\label{tab:comparison}
	\centering
	\begin{tabular}{c@{\hspace{4pt}}c@{\hspace{4pt}}c@{\hspace{4pt}}c}
		\hline\hline
		Method                                                    & \textbf{PM}    & \textbf{M1} & \textbf{M2}    \\
		\hline
		\tabincell{c}{Operating profit (10$^3$ CNY)}                       &   85.52        &  65.28      &  85.57          \\
		\tabincell{c}{Hydrogen yield (kg)}                        &   4303.0       &  3779.0     &  4304.3         \\
		\tabincell{c}{Grid electricity purchase (MWh)}                &   8.54         &  23.10      &  8.54          \\
		\tabincell{c}{Grid-code compliance}                      & \ding{51}       & \ding{51}  & \ding{55}    \\
		\tabincell{c}{Average 35-kV-side 11th/23rd \\harmonic currents (A)}   & 3.85/1.19        &4.18/0.98        &7.39/2.10     \\
		%  \tabincell{c}{Average 11th/23rd \\harmonic currents (A)}   & 13.46/4.18        &14.63/3.44        &25.88/7.36     \\
		\hline\hline
	\end{tabular}
\end{table}

\subsection{Large-Scale System Validation}

\subsubsection{Application to the 20-ELZ System}

\begin{figure}[t]
	\centering
	\includegraphics[width=3.45in]{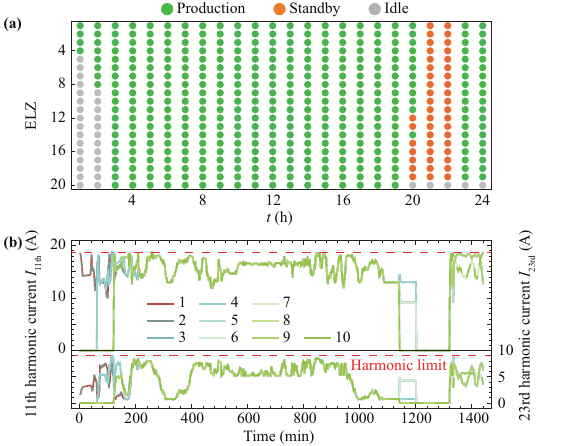}\vspace{-6pt}
	\caption{(a) Day-ahead on-off switching of ELZs and (b) real-time harmonic currents of each ELZ group in the 20-ELZ system.}
	\label{fig:large}
	\vspace{12pt}
	\includegraphics[width=3.45in]{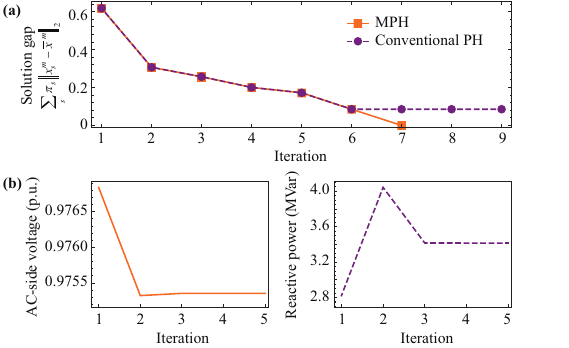}\vspace{-6pt}
	\caption{Convergence performance of the proposed method. (a) Day-ahead scheduling. (b) An illustrative real-time dispatch instance.}
	\label{fig:DHconvergence}
\end{figure}

To verify scalability, the PM is applied to the 20-ELZ system, with the ELZs divided into ten two-ELZ groups. Fig. \ref{fig:large} shows the day-ahead commitments and real-time harmonic currents. The upper layer still provides feasible ELZ commitments, while the lower-layer group-wise current and tap coordination keep the 11th and 23rd harmonics below their limits. These results demonstrate that the framework remains applicable in large-scale ReP2H systems.

\subsubsection{Computational Performance}
Fig. \ref{fig:DHconvergence} presents the convergence of the scheduling and real-time dispatch layers. % and Table \ref{tab:computational} compares the computational time in the 4-ELZ and 20-ELZ systems.
The MPH algorithm converges within a limited number of iterations and avoids the residual binary inconsistency observed with conventional PH.
Real-time HP-HM iterations also converge rapidly because the HM SPs are solved independently for each ELZ group and interact only through network voltage and reactive-power updates.
For the 4-ELZ and 20-ELZ systems, day-ahead scheduling requires 1.6 and 12 min, respectively, while the average real-time dispatch times are 4 and 20 s, well below the 2-min dispatch interval, confirming the computational feasibility of the proposed framework.

% The average solution time of each real-time dispatch instance is far shorter than the 2-min dispatch interval, confirming the computational feasibility of online operation.

%\begin{table}[t]
%\renewcommand{\arraystretch}{1.1}
%\centering
%\caption{Average Computational Time of the Proposed Method}\vspace{-8pt}
%\label{tab:computational}
%\begin{tabular}{ccc}
%\hline\hline
%System        & Day-ahead Scheduling & Real-time dispatch \\
%\hline
%4-ELZ system  & 1.6 min  & 4  s \\
%20-ELZ system & 12  min  & 20 s \\
%\hline\hline
%\end{tabular}
%\end{table}

\section{Discussion and Engineering Implications}
\label{sec:discussion}

\subsection{Engineering Implications for TR Selection}

\subsubsection{Difference Between 12- and 24-Pulse TRs}

The engineering tradeoff between 12-TR and 24-TR mainly involves harmonic performance, investment, losses, and DC-side ripple. A 24-TR provides stronger inherent harmonic cancellation but requires additional phase-shifting transformers and rectifier bridges, increasing capital cost and transformer losses. Its DC-side current ripple is also smaller, giving approximately 0.1\% higher rectification efficiency \cite{agredano2026experimental}. Nevertheless, studies have shown that the larger ripple of 12-TRs has negligible impact on ELZ degradation \cite{haleem2026impact}. Therefore, if harmonic limits can be met through coordinated operation, the main advantage of the 24-TR is its larger harmonic margin instead of hydrogen production.

\subsubsection{Economic Comparison}

The 12-TR and 24-TR schemes are compared in the 4-ELZ system through annual simulations using twelve typical days from \cite{zeng2025optimal}. The results are summarized in Table \ref{tab:12-24}. The investment costs of 5-MW 12-TR and 24-TR units are set to 0.6  $\times10^6$ and 1.0 $\times10^6$ CNY, respectively \cite{gao2024advanced, zeng2025optimal}, and annualized over a 20-year lifetime at an 8\% discount rate. The higher 24-TR cost mainly arises from its phase-shifting RCT and additional rectifier bridges. Each vacuum OLTC costs 1.5 $\times10^5$ CNY and has a lifetime of $3\times10^5$ tap-changing actions, corresponding to 0.5 CNY/action.

As shown in Table \ref{tab:12-24}, the two schemes produce nearly identical amounts of hydrogen. The slight efficiency advantage of the 24-TR is largely offset by losses in the additional phase-shifting transformer. Their main differences are therefore harmonic levels, rectifier investment, and OLTC actions.

The coordinated 12-TR scheme requires 28,105 tap-changing actions per year across all RCTs, compared with 8030 actions for the 24-TR scheme. This higher duty reflects the larger inherent harmonics of 12-TRs and is consistent with the harmonic-violation periods observed in Fig. \ref{fig:M2}. Nevertheless, the additional annual OLTC action cost is only $0.101\times10^5$ CNY, whereas the annualized rectifier investment decreases by $1.630\times10^5$ CNY. The resulting net saving is $1.529\times10^5$ CNY/year, corresponding to a $\frac{1.630-0.101}{4.074}=37.5\%$ reduction in rectification-stage cost relative to the 24-TR scheme.

The economic boundary is also favorable to the 12-TR. Under the typical OLTC lifetime and action cost, the 12-TR remains preferable as long as its investment is more than $2.5\times10^4$ CNY lower per unit than that of the 24-TR. This threshold is much smaller than the market cost difference. Thus, coordinated 12-TRs provide a practical alternative to higher-pulse rectification.

\begin{table}[t]
	\renewcommand{\arraystretch}{1.2}
	\caption{Performance Comparison of 12- and 24-TR in the 4-ELZ System}\vspace{-8pt}
	\label{tab:12-24}
	\centering
	\begin{tabular}{c@{\hspace{4pt}}c@{\hspace{4pt}}c@{\hspace{4pt}}c}
		\hline\hline
	     Scheme                                                   & \textbf{24-TR}    & \textbf{12-TR}   & Comparison    \\
		\hline
		\tabincell{c}{Operational revenue (10$^7$ CNY/yr)}           &   2.770        &  2.770     &  /          \\
		\tabincell{c}{Hydrogen yield (10$^6$ kg/yr)}                 &   1.337        &  1.337     &  /         \\
		\tabincell{c}{Annualized rectification-stage\\ CAPEX (10$^5$ CNY/yr)}                &   4.074        &  2.444     &  --1.630         \\
        \tabincell{c}{Action cost of RCT OLTCs (10$^5$ CNY/yr)}                &   0.040        &  0.141     &  +0.101         \\
        \tabincell{c}{Grid-code compliance}                      & \ding{51}       & \ding{51}  &      \\
		\tabincell{c}{Average 35-kV-side 11th/23rd \\harmonic currents (A)}   & 0/1.49       &  3.85/1.19        &     \\
%		\tabincell{c}{Average 11th/23rd \\harmonic currents (A)}   & 13.46/4.18        &14.63/3.44        &25.88/7.36     \\
		\hline\hline
	\end{tabular}
\end{table}

\subsection{Robustness to PCC SCC Variations}

\begin{figure}[t]
  \centering
  \includegraphics[width=3.5in]{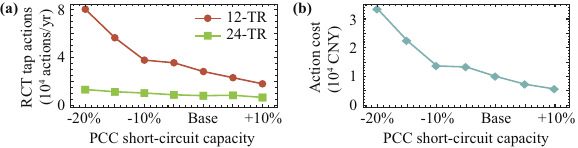}\vspace{-6pt}
  \caption{(a) OLTC action counts of 12/24-TR schemes and (b) additional OLTC action cost of 12-TR compared with 24-TR under different PCC SCCs.}
  \label{fig:sensitivity}
\end{figure}

The PCC short-circuit capacity determines the allowable harmonic-current injection and therefore affects the relative value of 12-TR and 24-TR schemes. A sensitivity analysis is conducted by varying the SCC from $-20\%$ to $+10\%$ of its base value. The resulting OLTC action counts and additional costs are shown in Fig. \ref{fig:sensitivity}.

Within this range, both schemes remain harmonic compliant through RCT OLTC adjustment without reducing hydrogen production. Their hydrogen yield and operating profit remain nearly unchanged, and the economic difference is dominated by rectifier investment and OLTC duty. As SCC decreases, tighter harmonic limits require more tap changes, particularly for the 12-TR because of its larger inherent harmonic injection. A stronger grid relaxes these limits and reduces the required OLTC actions, further favoring the 12-TRs.

\section{Conclusions}
\label{sec:conclusion}

This paper investigates harmonic-compliant operation of ReP2H systems using 12-pulse TRs and evaluates their potential as an alternative to higher-pulse TRs. A two-layer framework coordinates ELZ scheduling, renewable power tracking, and RCT OLTC-based harmonic mitigation. The main findings are:

\begin{enumerate}
	\item  RCT tap adjustment provides an additional degree of freedom for harmonic phasor regulation by modifying firing and commutation-overlap angles while preserving electrolytic-current references. Plant-level harmonic performance therefore depends on both rectifier characteristics and coordinated operating states.
	
	\item Joint current and tap regulation reduces the production penalty of harmonic compliance. HFR-based scheduling preserves sufficient real-time harmonic-regulation margin, while tap-assisted mitigation keeps currents close to their production-optimal allocation, reconciling renewable tracking, harmonic compliance, and voltage security.
	
	\item The 12-/24-pulse comparison reveals a practical tradeoff between higher-pulse hardware and operational harmonic flexibility. Under typical conditions, the cost of additional OLTC duty is lower than the avoided higher-pulse investment, and the tradeoff varies with grid strength.
\end{enumerate}

Future work will extend the framework to multi-stack plants with shared BoP systems and coupled hydrogen-side dynamics, including temperature and HTO impurity. Their interaction with electrical constraints will be incorporated into integrated scheduling and control to quantify harmonic-regulation flexibility under tighter process and safety constraints.

%\bibliographystyle{IEEEtran}
%\bibliography{IEEEabrv,12to24}

% Generated by IEEEtran.bst, version: 1.13 (2008/09/30)

\end{document}